\documentclass[12pt]{article}
\usepackage[cp866]{inputenc}
\usepackage{citehack}
\usepackage{amsmath}
\usepackage{amsfonts}
\usepackage{amssymb}
\usepackage{latexsym}
\newtheorem{theorem}{Theorem}

\def\zr{\ltimes}

\def\ZR{\rightthreetimes}

\def\spa{\text{\rm span}}

\def\Real{\mathbb{R}}
\def\g{\mathfrak{g}}
\def\h{\mathfrak{h}}
\def\so{\mathfrak{so}}

\def\z{\mathfrak{z}}

\def\Conf{\text{\rm Conf }}
\def\Simil{\text{\rm Sim }}
\def\Isom{\text{\rm Isom }}
\def\A{\mathcal {A}}
\def\B{\mathcal {B}}

\def\N{\mathcal {N}}
\def\K{\mathcal {K}}
\def\id{\text{\rm id }}

\title{Isometry groups of Lobachevskian spaces, similarity transformation
groups of Euclidean spaces and \\ Lorentzian holonomy groups}

\author{Anton S. Galaev}
\begin{document}

\maketitle\vskip-50ex
 {\renewcommand{\abstractname}{Abstract}
\begin{abstract}
Weakly-irreducible not irreducible subalgebras of $\so(1,n+1)$
were classified by \\ L. Berard Bergery and A. Ikemakhen.
In the present paper a geometrical
proof of this result is given.
Transitively acting isometry groups
of Lobachevskian spaces and transitively acting
similarity transformation groups of Euclidean
spaces are classified.

\end{abstract}

\section*{Introduction}

In 1952 A. Borel and A. Lichnerowicz showed that {\it the holonomy
group of a Riemannian manifold is a product of irreducible
holonomy groups of Riemannian manifolds}, see \cite{Bo-Li}. The main reason is the
following. If a subgroup $G\subset SO(n)$ preserves a proper
vector
subspace, then $G$ preserves also its orthogonal complement
$U^\bot$ and we have $R^n=U\oplus U^\bot$, i.e. the group $G$ is totally reducible.
In 1955 M. Berger classified  possible connected irreducible
holonomy groups of Riemannian manifolds, see \cite{Ber}.

The Borel and  Lichnerowicz theorem does not work in  the
pseudo-Riemannian case. Suppose a subgroup  $G\subset SO(r,s)$ preserves a
proper vector subspace $U\subset \Real^{r,s}$ such that the
restriction of the inner product to $U$ is degenerate, then
$U\cap U^\bot\neq\{0\}$ and we have no orthogonal decomposition of
$\Real^{r,s}$ into $G$-irreducible subspaces. A subgroup $G\subset SO(r,s)$
is called {\it weakly-irreducible} if it preserves no nondegenerate
proper subspace of $\Real^{r,s}$.
There is Wu's theorem that states that {\it the holonomy group of
a pseudo-Riemannian manifold is a product of
weakly-irreducible holonomy groups of pseudo-Riemannian manifolds},
see \cite{Wu}. If a
holonomy group is irreducible, then it is weakly-irreducible. In \cite{Ber}
M. Berger gave  a classification of irreducible holonomy groups for
pseudo-Riemannian manifolds. In particular, {\it the
only connected irreducible holonomy group of Lorentzian
manifolds is  $SO^0(1,n+1)$}, see \cite{Disc-Ol} and \cite{Bo-Ze} for
direct proofs of this fact.

There is still no classification of weakly-irreducible not
irreducible holonomy groups of pseudo-Riemannian manifolds. The
first step towards a classification of weakly-irreducible not
irreducible holonomy groups of Lorentzian manifolds was made by
L. Berard Bergery and A. Ikemakhen, who classified
weakly-irreducible not irreducible subalgebras of $\so(1,n+1)$,
see \cite{B-I}. 
More precisely, they divided weakly-irreducible not irreducible subalgebras of $\so(1,n+1)$
into 4 types.
The proof of this result  was purely algebraical.

We introduce a  geometrical proof of the result of L. Berard Bergery and A. Ikemakhen. 
We consider an $n+2$-dimensional Minkowski vector space $(V,\eta)$ and fix an isotropic vector $p\in V$.
We denote by $SO(V)_{\Real p}$
the Lie subgroup of $SO(V)$ that preserves the isotropic
line $\Real p$. We  denote by $E$ a vector subspace
$E\subset V$ such that $(\Real p)^{\bot_\eta}=\Real p\oplus E$, and
by $q$ an isotropic vector $q\in
V$ such that $\eta(q,E)=0$ and $\eta(p,q)=1$.
The vector space $E$ is an Euclidean space.
We consider the
vector model of the $n+1$-dimensional Lobachevskian space
$L^{n+1}$
and its boundary $\partial L^{n+1}$, which is diffeomorphic to the
$n$-dimensional unit sphere.  We have the natural isomorphisms
$$SO(V)\simeq\Isom L^{n+1}\simeq\Conf \partial  L^{n+1} \text{ and }
SO(V)_{\Real p}\simeq\Simil E,$$ where $\Isom L^{n+1}$
is the group of all isometries of $L^{n+1}$,
$\Conf\partial L^{n+1}$ is the group of all conformal transformations of
$\partial L^{n+1}$ and $\Simil E$ is the group of all similarity transformations of
$E$. We identify the set $\partial L^{n+1}\backslash\{\Real p\}$ with
the Euclidean space $E$. Then any subgroup $G\subset SO(V)_{\Real p}$
acts on $E$, moreover we have $G\subset\Simil E$.
We prove that {\it a connected  subgroup $G\subset SO(V)_{\Real p}$ is weakly-irreducible
iff the corresponding subgroup
$G\subset\Simil E$ under the isomorphism  $SO(V)_{\Real p}\simeq\Simil E$ acts transitively on $E$.}
This gives us a {\it one-to-one correspondence between connected
weakly-irreducibly acting subgroups of $SO(V)_{\Real p}$ and connected  transitively acting
subgroups of $\Simil E$.}
Using a description for connected transitive subgroups of 
$\Simil E$ (see \cite{Al2}, \cite{A-V-S}),
we divide such subgroups into 4 types.
We show that the corresponding weakly-irreducible subgroups of 
$SO(V)_{\Real p}$
have the same type introduced by L. Berard Bergery and A. Ikemakhen.

We also classify transitively acting isometry groups
of the Lobachevskian space $L^{n+1}$. We show that these
groups are exhausted by $SO^0(V)$ and by the weakly-irreducible not irreducible subgroups of
$SO(V)_{\Real p}$  of type 1 and type 3.

{\bf Remark} In another paper we will use a similar ideas
for complex Lobachevskian space in order to classify connected weakly-irreducible not irreducible
subgroups of $SU(1,n+1)\subset SO(2,2n+2)$.

{\bf Acknowledgments.} I wish to thank D.V. Alekseevsky for his
useful suggestions. Also I would like to express my gratitude to Helga Baum and M.V. Losik for
help in preparation of this paper.


\section{Results of L. Berard Bergery and A. Ikemakhen}

Let $(V,\eta)$ be a Minkowski space of dimension $n+2$,
where  $\eta$ is a metric on $V$ of
signature $(1,n+1)$. We fix a basis
$p,e_1,...,e_n,q$ of $V$ with respect to which the Gram matrix
of $\eta$ has the form
$\left(\begin{array}{ccc}
0 & 0 & 1\\ 0 & E_n & 0 \\ 1 & 0 & 0 \\
\end{array}\right)$,
 where $E_n$ is the $n$-dimensional identity matrix.

Let $E\subset V$ be the vector subspace spanned by $e_1,...,e_n$.
The vector space $E$ is an Euclidean space with respect to the inner product
$\eta|_E.$

Denote by $\so(V)$ the Lie algebra of all $\eta$-skew symmetric
endomorphisms of $V$ and by $\so(V)_{\Real p}$ the subalgebra
of $\so(V)$ that preserves the line $\Real p$.

The Lie algebra $\so(V)_{\Real p}$ can be identified with
the following algebra of matrices:
  $$\so(V)_{\Real p}=\left\{ \left (\begin{array}{ccc}
a &-X^t & 0\\ 0 & A & X \\ 0 & 0 & -a \\
\end{array}\right):\, a\in \Real,\, X\in E,\,A \in \so(E)
 \right\}.$$

The above matrix can be identified with the triple $(a,A,X)$.
Define the following subalgebras of $\so(V)_{\Real p}$,
 $\mathcal A=\{(a,0,0):a\in \Real\},$ $\mathcal
K=\{(0,A,0):A\in \mathfrak{so}(E)\}$ and
$\mathcal N=\{(0,0,X):X\in E\}$. 
We see that $\A$ commutes with $\K$, and $\N$ is an ideal.
We have the decomposition
$$\so(V)_{\Real p}=(\mathcal A\oplus\mathcal K)\zr\mathcal
N.$$

A subalgebra $\g\subset \so(V)$ is
called {\it irreducible} if it preserves no proper
subspace of $V$; $\g$ is called {\it weakly-irreducible} if
it preserves no nondegenerate proper subspace of $V$.

Obviously, if $\g\subset \so(V)$ is
irreducible, then it is weakly-irreducible.
If $\g\subset \so(V)$ preserves
a degenerate proper subspace $U\subset V$, then it preserves the
isotropic line $U\bigcap U^\bot$;
any such algebra is conjugated to a subalgebra of $\so(V)_{\Real p}$.

Let $\B\subset\so(E)$ be a subalgebra. Recall that $\B$ is a compact Lie
algebra and we have the decomposition $\B=\B'\oplus\z(\B)$, where
$\B'$ is the commutant of $\B$ and $\z(\B)$ is the center of $\B$.

The following result is due to L. Berard Bergery and A. Ikemakhen.

{\bf Theorem}
{\it Suppose $\g\subset \so(V)_{\Real p}$ is a weakly-irreducible subalgebra.
Then $\g$ belongs to one of the following types
\begin{description}
\item[type 1.] $\g=(\mathcal A\oplus\mathcal B) \zr\mathcal N$, where
$\B\subset\so(E)$ is a subalgebra;
\item[type 2.] $\g=\mathcal B \zr\mathcal N$;
\item[type 3.] $\g=(\mathcal B'\oplus \{\varphi(A)+A:A\in \mathcal \z(\B)\})\zr \N$, 
where $\varphi :\z(\B)\to\A$ is a non-zero linear map;
\item[type 4.] $\g=(\mathcal B'\oplus\{\psi(A)+A:A\in \z(B)\})\zr \N_W$, 
where we have a non-trivial   orthogonal
decomposition $E=U\oplus W$ such that $\B\subset\so(W)$;
$\N_W=\{(0,0,X):X\in W\}$;\\ $\N_U=\{(0,0,X):X\in U\}$ and
$\psi:\z(\B)\to \mathcal \N_U$ is a surjective linear map.
\end{description}}

Denote by $SO(V)$ the Lie group of all automorphisms of $V$ that
preserve the form $\eta$, and  with $\det f=1,$
and by $SO(V)_{\Real p}$ the  Lie subgroup of
$SO(V)$ that preserves the isotropic line $\Real p$. Obviously, $\so(V)$ and
$\so(V)_{\Real p}$ are the  Lie algebras of $SO(V)$
and $SO(V)_{\Real p}$ respectively.

By definition, the {\it type} of a connected weakly-irreducible
  Lie subgroup $G\subset SO(V)_{\Real p}$
is the type of its Lie algebra $\g\subset\so(V)_{\Real p}$.

\section{Transitive  similarity transformation groups of
Euclidean spaces}

In this section we recall a description for connected
transitively acting groups of similarity transformations and isometries of
Euclidean spaces, see \cite{Al2} or \cite{A-V-S}.

Let $(E,\eta)$ be an Euclidean space. A map $f:E\to E$
is called a {\it similarity transformation} of $E$ if there exists
a $\lambda >0$ such that $\|f(x_1)-f(x_2)\|=\lambda\|x_1-x_2\|$ for all
$x_1,x_2\in E$, where $\|x\|^2=\eta(x,x)$. If $\lambda=1$, then $f$
is called an {\it isometry}.
Denote by $\Simil E$ and $\Isom E$ the groups of all similarity
transformations and isometries of $E$ respectively.
A subgroup $G\subset \Simil E$ such that $G\not\subset\Isom E$ is called
{\it essential}. A subgroup $G\subset \Simil E$ is called {\it irreducible}
if it preserves no proper affine subspace of $E$.

We need the following theorem from \cite{Al2} and \cite{A-V-S}.
 
\begin{theorem}\label{AVS}
{\rm (1)} Let $G\subset \Isom E$ be a connected subgroup that
acts transitively on $E$. 
Then there exists a
decomposition $G=H\ZR F$, where $H$ is the stabilizer of a point
$x\in E$ and $F$ is a normal subgroup of $G$ that acts simply
transitively on $E$.

{\rm (2)} Let $F\subset \Isom E$ be a connected subgroup that acts
simply transitively on $E$. Then there exists an orthogonal decomposition
$E=U\oplus W$ and a Lie groups homomorphism\\  $\Psi:U\to SO(W)$
such that $F=U^\Psi\ZR W$, where 
$$U^\Psi=\{\Psi(u)\cdot u:u\in U\}\subset SO(W)\times U$$ is a
group of screw isometries.

{\rm (3)} 
Let $G\subset\Simil E$ be an essential connected subgroup that  acts
 transitively on $E$. Then $G=(A_1\times H)\ZR F$,
where $A_1\subset \Simil E$ is a 1-parameterized essential
subgroup that preserves a point
$x$, $H\subset\Isom E$ commutes with $A_1$ and preserves the point $x$, 
and $F$ is a normal
subgroup of $G$ that acts simply transitively on $E$.

{\rm (4)} A connected subgroup $G\subset\Isom  E$ acts irreducibly on $E$ if
and only if it acts transitively on $E$.
\end{theorem}

From parts (3) and (4) of the theorem it follows that 
{\it a connected  subgroup $G\subset\Simil  E$ acts irreducibly on $E$ if
and only if it acts transitively on $E$.}

\section{Isometries of Lobachevskian spaces}

Let $p,e_1,...,e_n,q$ be a basis of the vector space $V$ as
above.  
Consider the basis $e_0,e_1,...,e_n,e_{n+1}$ of $V$, 
where  $e_0=\frac{\sqrt{2}}{2}(p-q)$ and $e_{n+1}=\frac{\sqrt{2}}{2}(p+q)$.
With respect to this basis the Gram matrix of $\eta$ has the form
 $\left(
\begin{array}{cc}
-1 & 0 \\
0 & E_{n+1}\\
 \end{array} \right),$ where $E_{n+1}$ is the $n+1$-dimensional
 identity matrix.

The vector model of the $n+1$-dimensional {\it Lobachevskian space} is
defined  in the following way
$$L^{n+1}=\{x\in V:\ \eta(x,x)=-1,\ x_0>0\}.$$
Recall that $L^{n+1}$ is an $n+1$-dimensional Riemannian submanifold of $V$.
The tangent space at a point $x\in L^{n+1}$ is identified with the
vector subspace $(x)^{\bot_\eta}\subset V$ and the restriction of the form
$\eta$ to this subspace is positively definite.

Any element $f\in SO(V)$ preserves the space $L^{n+1}$.
Moreover, for any $f\in SO(V)$, the restriction $f|_{L^{n+1}}$ is
an isometry of $L^{n+1}$ and any isometry of $L^{n+1}$ can be
obtained in this way. Hence we have the isomorphism
$$SO(V)\simeq\Isom L^{n+1},$$
where $\Isom L^{n+1}$ is the group of all isometries of
$L^{n+1}$.

Consider the {\it light-cone} of $V$,
$$C=\{x\in
V:\ \eta(x,x)=0\}.$$
The subset of the $n+1$-dimensional projective space $\mathbb{P}V$ that
consists of all {\it isotropic lines} $l\subset C$ is called the 
{\it boundary of the Lobachevskian space} $L^{n+1}$ and is 
denoted by $\partial L^{n+1}$.

We identify $\partial L^{n+1}$ with the $n$-dimensional unit sphere $S^n$ in
the following way.
Consider the vector subspace $E_1=E\oplus\Real e_{n+1}$. 
Each isotropic line intersects the hyperplane $e_0+E_1$ at a unique
point. The
intersection $(e_0+E_1)\cap C$ is the set $$\{x\in V : x_0=1,\,
x_1^2+\cdots+ x_{n+1}^2=1\},$$ which is the $n$-dimensional sphere $S^n$.
This gives us the identification $\partial L^{n+1}\simeq S^n$.

Denote by $\Conf S^n$ the group of all conformal transformations
of $S^n$. Any transformation $f\in SO(V)$ takes isotropic lines to
isotropic lines. Moreover, under the above identification, we have
$f|_{\partial L^{n+1}}\in \Conf \partial L^{n+1}$ and any transformation from
$\Conf \partial L^{n+1}$ can be obtained in this way. Hence we have the
isomorphism $$SO(V)\simeq\Conf\partial L^{n+1}.$$

Suppose $f\in SO(V)_{\Real p}$. The corresponding element
$f\in\Conf S^n$ (we denote it by the same letter) preserves the
point $p_0=\Real p\cap (e_0+ E_1)$. Clearly, $p_0=\sqrt 2 p$.
Denote by $s_0$ the stereographic projection $s_0:
S^n\backslash\{p_0\}\to e_0+E$.
Since $f\in\Conf S^n$, we see that $s_0\circ f \circ s^{-1}_0:E\to E$
(here we identify $e_0+E$ with $E$) is a similarity transformation of
the Euclidean space $E$. Conversely, any similarity transformation of
$E$ can be obtained in this way. Thus we have the isomorphism
$$SO(V)_{\Real p}\simeq \Simil E.$$

A {\it plane} in the Lobachevskian space $L^{n+1}$ is the nonempty
intersection of $L^{n+1}$ and of a vector subspace $U\subset V$.
The intersection $L^{n+1}\cap U$ is not empty if and only if
the restriction of the form $\eta$ to $U$ has   signature $(1,\dim U-1)$. 
A subgroup $G\subset \Isom L^{n+1}$ is called {\it irreducible} if
it preserves no proper plane in $L^{n+1}$.

The following theorem is due to F.I. Karpelevich, see \cite{A-V-S} or
\cite{Kar}.
\begin{theorem} \label{Kar}  Let $G$ be a proper connected closed subgroup
of $\Isom L^{n+1}$. Then $G$ acts irreducibly on $L^{n+1}$ if and
only if it preserves an isotropic line $l\in\partial L^{n+1}$ and
acts transitively on the Euclidean space 
$E_l=\partial L^{n+1}\backslash\{l\}$.
\end{theorem}

Since the holonomy group of a Lorentzian manifold can be
not closed, we need an analog of this theorem for not
closed groups.
In \cite{Disc-Ol} was proved the following theorem.
\begin{theorem} \label{DiScOl} 
Let $G$ be a connected (non nec. closed) subgroup
of $SO(V)$ that acts weakly-irreducibly. 
Then either $G$ acts transitively on $L^{n+1}$ or
$G$ acts transitively on the Euclidean space 
$E_l=\partial L^{n+1}\backslash\{l\}$.
\end{theorem}
We will prove the following theorem.
\begin{theorem} \label{wir} Let $G$ be a proper connected  subgroup
of $SO(V)_{\Real p}$. Then $G$ acts weakly-irreducibly on $V$ if and
only if it acts transitively on the Euclidean space 
$E=\partial L^{n+1}\backslash\{\Real p\}$.
\end{theorem}
{\bf Proof.} We claim that the subgroup
$G\subset SO(V)_{\Real p}$ acts weakly-irreducibly on $V$ if and
only if the corresponding subgroup $G\subset\Simil E$ acts irreducibly on $E$.
If $G\subset SO(V)_{\Real p}$ is not  weakly-irreducible, then
it preserves a not degenerate proper subspace $U\subset V$.
Since the orthogonal complement $U^\bot\subset V$ is also preserved 
and either $U\cap C\neq\{0\}$ or  $U^\bot\cap C\neq\{0\}$,
we can assume that  $U\cap C\neq\{0\}$.
The subgroup $G\subset\Simil E$ preserves the affine
subspace $s_0((e_0+E)\cap C\cap U)\subset E$, which is not empty.
Conversely, if the subgroup $G\subset\Simil E$ preserves a proper affine
subspace $W\subset E$, then  $G\subset SO(V)_{\Real p}$ preserves
the vector subspace of $V$ spanned by $s_0^{-1}(W)\subset e_0+E$, which is not degenerate.
Now the proof of the theorem follows from
parts (3) and (4) of theorem \ref{AVS}. $\Box$

\section{Application to holonomy groups of Lorentzian manifolds}

Now we consider connected  weakly-irreducible not irreducible subgroups of $SO(V)$.
Any such group $G$ preserves an isotropic line and is conjugated
to a subgroup of $SO(V)_{\Real p}$.

In section 2 we have constructed the isomorphism
$SO(V)_{\Real p}\simeq \Simil E$.
This isomorphism and theorem \ref{wir}  gives us a {\it one-to-one correspondence between connected
weakly-irreducible subgroups $G\subset SO(V)_{\Real p}$ and connected transitively acting
subgroups $G\subset\Simil E$.}
\begin{theorem} \label{hol}  Let $G\subset\Simil E$ be a
transitively acting connected subgroup. Then $G$ belongs to one of the following types
\begin{description}
\item[type 1.] $G=(A\times H)\ZR E$, where $A=\Real^+$ is the unite component
for the group of all dilations of $E$ about the origin $0$,
$H\subset SO(E)$ is a Lie subgroup, and $E$ is the group of all translations in $E$;
\item[type 2.] $G=H\ZR E$;
\item[type 3.] $G=(A^\Phi \times H)\ZR E$, 
where $\Phi: A\to SO(E)$ is a homomorphism and 
$$A^\Phi=\{\Phi(a)\cdot a:a\in A\}\subset SO(E)\times A$$ 
 is a group of screw
dilations of $E$;
\item[type 4.] $G=(H\times U^{\Psi})\ZR W,$ where $E=U\oplus W$ is an orthogonal decomposition,
$\Psi:U\to SO(W)$ is a homomorphism, and
$$U^{\Psi}=\{\Psi(u)\cdot u:u\in U\}\subset SO(E)\times U$$
is a group of screw isometries of $E$.
\end{description}

The corresponding
subgroups of $SO(V)_{\Real p}$ under the isomorphism
$SO(V)_{\Real p}\simeq \Simil E$ are the groups
of the same type introduced by L. Berard Bergery and A. Ikemakhen.
\end{theorem}

{\bf Proof.} Denote by $A$, $K$ and $N$ the connected Lie subgroups of
$SO(V)_{\Real p}$ corresponding to the subalgebras $\A$, $\K$
and $\N\subset\so(V)_{\Real p}$.
With respect to the basis $p,e_1,...,e_n,q$ 
these groups have the following forms
$A=\left\{ \left (\begin{array}{ccc}
a & 0 & 0\\ 0 & \id & 0\\ 0 & 0 & \frac{1}{a}\\
\end{array} \right):a\in \Real,\,a>0 \right\},$

$K=\left\{ \left(
\begin{array}{ccc}
1 & 0 & 0\\ 0 & f & 0\\ 0 & 0 & 1 \\
\end{array} \right):\,f\in SO(E) \right \}$ and
$N=\left\{ \left ( \begin{array}{ccc} 
1 &-X^t & -\frac{1}{2}X^tX\\ 0 & \id & X\\ 0 & 0& 1 \\
\end{array} \right):X\in E \right \}.$

We have the decomposition $SO^0(V)_{\Real p}=(A\times K)\ZR N$.

The computation shows that under the isomorphism
$SO(V)_{\Real p}\simeq \Simil E$\\
the element
$\left (\begin{array}{ccc}
a & 0 & 0\\ 0 & \id & 0\\ 0 & 0 & \frac{1}{a} \\
\end{array} \right)\in A$
\begin{tabular}{ll}
corresponds to the dilation & $X\mapsto aX$,\\
\end{tabular}

the element
$\left (
\begin{array}{ccc}
1 & 0 & 0\\ 0 & f & 0\\ 0 & 0 & 1 \\
\end{array} \right)\in K$
corresponds to $f\in SO(E)$, and

the element
$\left ( \begin{array}{ccc}
 1 & -X^t & -\frac{1}{2}X^tX\\ 0 & \id & X\\ 0 & 0& 1 \\
\end{array} \right)\in N$ 
\begin{tabular}{ll}
corresponds to the translation & $Y\mapsto Y+X$.\\
\end{tabular}

Let a subgroup $G\subset\Simil E$ act transitively. 
Denote by the same letter $G$ the corresponding
weakly-irreducible subgroup of $SO(V)_{\Real p}$.
Since we are interested in the groups up to
conjugacy, in the theorem \ref{AVS} we choose $x=0$, 
then $H\subset SO(E)$.

For the subgroup $G\subset SO(V)_{\Real p}$ we have two cases: 

{\bf case 1.} $G$
preserves the vector $p$; 

{\bf case 2.} $G$ preserves the isotropic line
$\Real p$ but does not preserve the vector $p$.

Consider these cases.

{\bf Case 1.} We have $G\subset K\ZR N$. Hence the corresponding subgroup
$G\subset\Simil E$ consists of isometries, i.e. $G\subset \Isom E$.
From the transitivity of $G$ it follows that $G=H\ZR F$, where
$H\subset SO(E)$ and $F$ is a normal subgroup of $G$  that acts simply
transitively on $E$. Hence there exists an orthogonal decomposition
$E=U\oplus W$ and a homomorphism
$\Psi:U\to SO(W)$ such that $F=U^{\Psi}\ZR W$. 

There are two subcases

{\bf Subcase 1.1.} The homomorphism $\Psi$ is trivial. Hence $F=E$
and $G=H\ZR E$. From the classification of L. Berard Bergery and A. Ikemakhen
we have $G\subset SO(V)_{\Real p}$ is {\it a group of type 2}.

{\bf Subcase 1.2.} The homomorphism $\Psi$ is not trivial.
We can assume that the homomorphism $d\Psi:U\to\so(W)$
is injective. Indeed, if $\ker d\Psi\neq\{0\}$, then we choose the
decomposition $E=U_1\oplus W_1$, where $W_1=W\oplus\ker d\Psi$ and
$U_1\subset U$ is the orthogonal complement of $\ker d\Psi$ in
$U$, and we consider $\Psi_1=\Psi|_{U_1}$.

We claim that $H$ commutes with $\Psi(U)\subset SO(W)$, moreover
$H$ acts trivially on $U$ and $H\subset SO(W)$.
Let $f\in H$, $u\in U$. Since $F$ is a normal subgroup of $G$,
we have $f\circ \Psi(u)\circ u\circ f^{-1}=w\circ \Psi(u_1)\circ u_1$ 
for some $w\in W$ and $u_1\in U$.
Hence  for all $v\in E$ we have
$f(u)+f\circ \Psi(u)\circ f^{-1}(v)=w+u_1+\Psi(u_1)v$.
Since this holds for all $v\in E$, we have  
$f\circ \Psi(u)\circ f^{-1}=\Psi(u_1)$.
We will prove that $\Psi(u)=\Psi(u_1)$.
Let $l(\Psi(U))$ and $\h=l(H)$ be the Lie algebras of the Lie groups
$\Psi(U)$ and $H$ respectively.
We have $(\h+l(\Psi(U)))'=\h'+[\h,\Psi(U)]$.
Since $[\h,\Psi(U)]\subset \Psi(U)$ and the Lie algebra
$l(\Psi(U))$ is commutative, we have $(\h+l(\Psi(U)))''=\h'$.
If $\Psi(u)\neq\Psi(u_1)$, then $[\h,\Psi(U)]\neq\{0\}$ and
$(\h+l(\Psi(U)))'\neq(\h+l(\Psi(U)))''$. Since the 
subalgebra $\h+l(\Psi(U))\subset\so(E)$ is compact, we have a contradiction.
Thus, $\Psi(u)=\Psi(u_1)$ and $H$ commutes with $\Psi(U)$.
Consider now the Lie algebra $l(G)$ of the Lie group $G$.
We have $l(G)=(\h\oplus l(U^\Psi))\zr W$.
Since $U^{\Psi}=\{\Psi(u)\circ u:u\in U\}$, we see that
$l(U^{\Psi})=\{d\Psi(u)+u:u\in U\}$.
For $\xi\in \h$ and $d\Psi(u)+u\in l(U^{\Psi})$ 
we have $[\xi,d\Psi(u)+u]=\xi u\subset U$. Since $U\cap l(G)=\{\varnothing\}$,
we see that $\xi u=0$. Hence $H$ acts trivially on $U$. Since $H\subset SO(E)$
and $W$ is orthogonal to $U$, we see that $H(W)\subset W$ and $H\subset SO(W)$.

We see now that $d\Psi(U)\subset\so(W)$ is a commutative subalgebra
that  commutes with $\h$.
Put $\B=\h\oplus d\Psi(U)$. We have $\z(\B)=\z(\h)\oplus d\Psi(U)$.
Put $\psi=d\Psi^{-1}:d\Psi(U)\to U$ and extend $\psi$ to the linear
map $\psi:\z(\B)\to U$ by putting $\psi|_{\z(\h)}=0$.
Thus we have $$l(G)=(\B'\oplus\{\psi(A)+A:A\in \z(\B)\})\zr W.$$
We see that $l(G)$ is an algebra of type 4 and $G$ is {\it a group of
type 4}.

{\bf Case 2.} In this case we have $G\subset \Simil E$, hence $G=(A_1\times H)\ZR F$, where
$A_1$ is a 1-parameterized subgroup of $G$ that preserves the point
$0$, $H\subset SO(E)$ commutes with $A_1$, and $F$ is a normal
subgroup that acts simply transitively on $E$.

There are two subcases

{\bf Subcase 2.1.} We have $A_1=A$ is the unity component of the group
of all dilations of $E$ about the origin $0\in E$.

We claim that $F=E$. Indeed, suppose that $F=U^\Psi\ZR  W$ and the
homomorphism $\Psi$ is not trivial. Let $u\in U$, $w\in W$ and
$1\neq\lambda\in A=\Real^+$. Since the subgroup $F\subset G$ is  normal,
we see that $\lambda\circ
\Psi(u)\circ u\circ w\circ\lambda^{-1}\in U^{\Psi}\ZR W$.
Let $v\in E$. We have
$(\lambda\circ\Psi(u)\circ u\circ w\circ\lambda^{-1})v=$\\ $\Psi(u)(\lambda\circ
u\circ w\circ\lambda^{-1})v=\Psi(u)(\lambda\circ u\circ w(\lambda^{- 1}v))=\Psi(u)(
\lambda(u+w+\lambda^{-1}v))=\Psi(u)(\lambda u +\lambda w+v)$.\\
Hence, $\lambda\circ
\Psi(u)\circ u\circ w\circ\lambda^{-1}=\Psi(u)\circ (\lambda u)\circ (\lambda w)\in U^{\Psi}\ZR W$.
This implies $u=\lambda u$ for all $u\in U$, hence,
$\lambda=1$. This gives us a contradiction. Thus, $F=E$.

Now we see that $G=(A_1\times H)\ZR F$ is {\it a group of type 1.}

{\bf Subcase 2.2.} In this case $A_1\neq A$, then $A_1\subset A\times SO(E)$.
By analogy with subcase 2.1., we can prove that $F=E$.

Let $\xi:\Real\to A_1$ be a parameterization of the group $A_1$.
Define the homomorphisms $\xi_1:\Real \to A$ and $\xi_2:\Real \to SO(E)$
by condition $\xi(t)=\xi_1(t)\cdot\xi_2(t)$ for all $t\in\Real$.
Since $A_1\not\subset SO(E)$, we see that the homomorphism
$\xi_1$ is an isomorphism. Put
$\Phi=\xi_2\circ\xi_1^{-1}:A\to SO(E)$. We have
$$A_1=\{\Phi(a)\cdot a:a\in A\}\subset SO(n)\times \Real.$$

We see that $l(G)=(l(A_1)\oplus\h)\zr E$ and
$$l(A_1)=\{d\Phi(a)+a:a\in l(A)\}.$$
Note that the subalgebra $l(d\Phi(l(A)))\subset\so(E)$ is
commutative and commutes with $\h$.
Put $\B=\h\oplus l(d\Phi(l(A)))$. We see that $\z(\B)=\z(\h)\oplus
l(d\Phi(l(A))).$
Put $\varphi=(d\Phi)^{-1}:d\Phi(l(A))\to l(A)$ and extend
$\varphi$ to the linear map $\varphi:\z(\B)\to l(A)$ by putting
$\varphi|_{\z(\h)}=0$.
Thus, $$l(G)=(\B'\oplus\{\varphi(A)+A:A\in \z(\B)\})\zr E.$$
We see that $G$ is {\it a group of type 3}.
This completes the proof of the theorem. $\Box$.

\section{Transitive isometry  groups  of 
the Lobachevskian space $L^{n+1}$}

Recall that  we consider a   Minkowski space
$(V,\eta)$  of dimension $n+2$
and  a basis
$p,e_1,...,e_n,q$ of $V$ with respect to which the Gram matrix
of $\eta$ has the form
$\left (\begin{array}{ccc}
0 & 0 & 1\\ 0 & E_n & 0 \\ 1 & 0 & 0 \\
\end{array}\right)$,
 where $E_n$ is the $n$-dimensional identity matrix.
We consider the vector subspace $E\subset V$ spanned by $e_1,...,e_n$ as
an Euclidean space with respect to the inner product $\eta|_E.$
We denote by  $SO(V)_{\Real p}$ the subgroup
of $SO(V)$ that preserves the line $\Real p$.
For the Lie group $SO^0(V)_{\Real p}$ 
we have the decomposition  $SO^0(V)_{\Real p}=(A\times K)\ZR N$, where
with respect to the basis $p,e_1,...,e_n,q$ 
the groups $A$, $K$ and $N$ have the following matrix forms
$A=\left\{ \left (\begin{array}{ccc}
a & 0 & 0\\ 0 & \id & 0\\ 0 & 0 & \frac{1}{a} \\
\end{array} \right):a\in \Real,\,a>0 \right \},$\\
 $K=\left\{ \left (
\begin{array}{ccc}
1 & 0 & 0\\ 0 & f & 0\\ 0 & 0 & 1 \\
\end{array} \right):\,f\in SO(E) \right \}$ and
$N=\left\{ \left ( \begin{array}{ccc} 
1 &-X^t & -\frac{1}{2}X^tX\\ 0 & \id & X\\ 0 & 0& 1 \\
\end{array} \right):X\in E \right \}.$
\begin{theorem}\label{ltrans}
Let $G\subset SO(V)$ be a connected subgroup that acts
transitively on the Lobachevskian space $L^{n+1}$.
Then either $G=SO^0(V)$ or $G$ preserves an
isotropic line $l\subset V$ and there exists a basis
$p,e_1,...,e_n,q$ of $V$ as above such that $l=\Real p$ and  $G$ is one of the 
following groups 

{\rm (1)} $(A\times H)\ZR N$, where $H\subset K$ is a subgroup;

{\rm  (2)} $(A^\Phi\times  H)\ZR N$,  
 where $\Phi: A\to K$ is a not trivial homomorphism and
 $$A^\Phi=\{\Phi(a)\cdot a:a\in A\}\subset K\times A.$$

Moreover the groups of the form $A\ZR N$ and $A^\Phi\ZR  N$
exhaust all connected subgroups of $SO(V)$ that act
simply transitively on $L^{n+1}$.
\end{theorem}  

Note that $A$ is the group of translations in
$L^{n+1}$ along the line $h=(\Real p\oplus\Real q)\cap L^{n+1}$,
$K$ is the group of rotations about $h$, $N$ is the group of parabolic
translations along 2-dimension planes in  $L^{n+1}$ 
and $A^\Phi$ is a group of screw translations along the line $h$.

{\bf Proof.}
Suppose a subgroup $G\subset SO(V)$  acts
transitively on $L^{n+1}$.  
Then it preserves no plane in $L^{n+1}$, hence $G$ acts
weakly-irreducibly on $SO(V)$.
If $G$ acts irreducibly on $V$, then $G=SO^0(V)$, see \cite{Disc-Ol} or \cite{Bo-Ze}.

If $G$ acts weakly-irreducibly not irreducibly  on $V$,
then $G$ preserves an isotropic line $l\subset V$, we assume that
$l=\Real p$. Then $G$ is the group of type 1,2,3 or 4.

We claim that the subgroup $K\ZR N\subset SO(V)$ does not act
transitively on $L^{n+1}$. Indeed, any element of
$K\ZR N$ takes the vector $\frac{1}{2}p-q\in L^{n+1}$ to some vector
$u-q$, where $u\in \spa\{p,e_1,...,e_n\}$, hence there is no element of $K\ZR N$
that takes $\frac{1}{2}p-q\in L^{n+1}$ to $p-\frac{1}{2}q\in L^{n+1}$.
Hence the groups of type 2 and 4 does not act transitively on $L^{n+1}$.

We must prove that groups of type 1 and 3, i.e.
groups of the form $A\times H\ZR N$ and 
$A^\Phi\times  H\ZR N$
act  transitively on $L^{n+1}$.
Let $v=xp+\alpha+yq\in L^{n+1}$ and $w=xp+\beta+yq\in L^{n+1}$,
where $\alpha,\beta\in E$. Then we have $2xy+\eta(\alpha,\alpha)=-1$
and $2xy+\eta(\beta,\beta)=-1$. Let $X=\frac{\alpha-\beta}{y}$.
The element $\left ( \begin{array}{ccc} 
1 & -X^t & -\frac{1}{2}X^tX\\ 0 & \id & X\\ 0 & 0& 1 \\
\end{array} \right)\in N$ takes $u$ to $w$.

Let $v=x_1p+\beta+y_1q\in L^{n+1}$, i.e.   $2x_1y_1+\eta(\beta,\beta)=-1$.

The element $\left (\begin{array}{ccc}
\frac{x_1}{x} & 0 & 0\\ 0 & \id & 0\\ 0 & 0 & 
\frac{x}{x_1} \\
\end{array} \right)\in A$ takes $w$ to $v$.
 The element $\left (\begin{array}{ccc}
\frac{x_1}{x} & 0 & 0\\ 0 & \Phi(\frac{x_1}{x}) & 0\\ 0 & 0 & 
\frac{x}{x_1} \\
\end{array} \right)\in A^\Phi$ takes $w$ 
to $xp+\Phi(\frac{x_1}{x})(\beta)+yq\in L^{n+1}$.
Thus there exist elements in $(A\times H)\ZR N$ and $(A^\Phi\times  H)\ZR N$ that take 
$u$ to $v$, i.e. the groups $(A\times H)\ZR N$ and  $(A^\Phi\times  H)\ZR N$ act
transitively on $L^{n+1}$.

Note that the elements of the subgroup 
$H\subset G$ preserve the point $p-\frac{1}{2}q\in L^{n+1}$.
Since $\dim L^{n+1}=\dim (A\ZR N)=\dim (A^\Phi\ZR N)$ and $L^{n+1}$ is simply connected,
we see that the groups of the  form  $A\ZR N$ and $A^\Phi\ZR  N$ are the only
connected subgroups of $SO(V)$ that act simply transitively on $L^{n+1}$. 
$\Box$ 

\bibliographystyle{unsrt}

\end{document}